\documentclass[10pt]{article}

\usepackage[latin1]{inputenc}
\usepackage[T1]{fontenc}
\usepackage[english]{babel}
\usepackage{amsfonts}
\usepackage{amssymb}
\usepackage{amsmath}
\usepackage{amsthm}
\usepackage{graphicx,xcolor}

\newtheorem*{lemma}{Lemma}
\newtheorem{teorema}{Theorem}
\newtheorem{corollario}{Corollary}
\newtheorem*{t*}{Theorem}

\def\be{\begin{equation}}
\def\ee{\end{equation}}
\def\bea{\begin{eqnarray}}
\def\eea{\end{eqnarray}}
\def\ni{\noindent}
\def\nn{\nonumber}

\def\s{\sigma}
\def\t{\tau}
\def\E{\mathbb{E}}
\def\max{\operatorname{max}}

\newcommand{\OO}[1]{O \left(\frac{1}{#1}\right)}

\title{Universality in bipartite mean field spin glasses}
\author{Giuseppe Genovese \footnote{email: giuseppe.genovese@mat.uniroma1.it}\\  \\ \textit{Dipartimento di Matematica, Sapienza Universit\`a di Roma}\\
\textit{Piazzale Aldo Moro, 2, 00185 Roma, Italia}}

\date{\today}

\begin{document}
\maketitle

\begin{abstract}
\ni In this work we give a proof of universality with respect to the choice of the statistical distribution of the quenched noise, for mean field bipartite spin glasses. We use mainly techniques of spin glasses theory, as Guerra's interpolation and the cavity approach. As a direct conseguence of our results, we have a proof of $L_\infty$ convergence of the free energy of the Hopfield Model to its expectation value.
\end{abstract}

\section*{Introduction}

Although spin glass models, expecially the well known Sherrington Kirkpatrick (SK) model, have been largely investigated in the past years \cite{MPV}\cite{guerrasg}\cite{talabook}, limited attention has been payed for bipartite systems. Despite that, bipartite spin glasses have a very interesting matematical structure, quite similar in several aspects to the Hopfield Model for neural networks \cite{hop}\cite{bg}\cite{BGG1}, and find many applications in modelling (see \cite{BGG2} and references therein).

\ni The problem of universality is an important aspect of the theory of random matrices \cite{AGZ}\cite{Tbook}, and very recently new general results have been achieved \cite{TV}. From the spin glass perspective, universalty is usually considered not of primarly importance, since it is believed that the main characteristics of spin glasses are independent on the choice of the distribution of quenched noise. It is infact the case of the SK Model. On the other hand, Guerra's interpolation, that is the key ingredient of the proof of the Parisi formula by Talagrand \cite{RSB}\cite{TP}, works only with gaussian random interaction (it is based on the Wick rule). Therefore, the problem to justify the particular choice of gaussian interaction in SK model was dealt and solved \cite{talbern}\cite{CH}\cite{ch}.

\ni The same problem arises in dealing with bipartite spin glasses: a systematic mathematical study of such a model, by interpolation method, has been started in \cite{BGG2}, motivated by a well known analogy with the Hopfield Model \cite{talabook}\cite{Bovbook}\cite{bg}\cite{BGG1}, that is infact a special bipartite spin glass.

\ni It will be clear from what follows that, even though the mathematical structure of bipartite spin glasses and the Hopfield Model is the same for many aspects, their universality properties are quite different. In particular, although we have no complete proof of that, we believe that universality does not hold for the Hopfield Model, at least not in the whole phase diagram.

\ni The work is organized as follows:

\ni In Section 1 the general structure of bipartite models is given, even in the case of systems of (bounded) soft spin, although in the whole paper only dichotomic variables are used (the generalization is straightforward).

\ni In Section 2 we will state the results: the pressure of the bipartite spin glass model, with certain conditions on the noise, is close in distribuition to the one with gaussian noise, when the size of the system grows to infinity. Then we state that their are close also in $L_p$ norm for some $p$ to be specified later.

\ni Proofs of the statements are given in Section 3.

\ni Lastly, in Section 4 we will point out the results that can be extended to the Hopfield Model, in particular we give a proof of the $L_\infty$ convergence of the pressure to the quenched one. Furthermore, we will discuss the points where our strategy fails for such a model, trying to give some explanations.

\section{The Structure of Bipartite Models}

We will deal with a set of N i.i.d. random spin variables $\sigma_i$, $i=1,..., N$, with any probability distribution $\mu(\sigma)$, symmetric with compact support $[-L, L]$. In particular we notice that $\mathbb{E}_{\sigma}[\sigma]=0$, and for a given bounded function of spin $f(\sigma)$, we must have $\mathbb{E}_{\sigma}[f(\sigma)]\leq L[\sup_{\sigma\in[-L, L]} f(\sigma)]$.

\ni Let us consider now for every $N$ another set of i.i.d. random spin variables $\tau_{\mu}$, $\mu=1,...,K$, with any probability distribution $\nu(\tau)$ with the above properties, but in general $\mu(\sigma)$ and $\nu(\tau)$ may be different. Therefore we have two distinct sets (or parties hereafter) of different spin variables, and we let them interact via the hamiltonian:
\begin{equation}\label{eq:H1}
H_{N,K}(\xi; \sigma, \tau)=-\sqrt{\frac{2}{N+K}}\sum_{i=1}^{N}\sum_{j=1}^{K} \xi_{ij}\sigma_i\tau_j,
\end{equation}
where the $\xi_{i\mu}$ are also i.i.d. r.v., with $\mathbb{E}[\xi]=0$ and $\mathbb{E}[\xi^2]=1$,  \textit{i.e.} the quenched noise ruling the mutual interactions between parties. It is then defined a mean field bipartite spin glass model \cite{BGG2}. 

\ni For sake of simplicity, in what follows we deal with both parties formed by dichotomic variables (as usual, sums denote not normalized expectations).

\ni The partion function, the pressure and the free energy of the model are defined as
\begin{eqnarray}
Z_{N, K}(\beta, \xi)&=&\sum_{\sigma,\tau}\exp\left(-\beta H_{N,K}(\xi; \sigma, \tau)\right),\\
A_{N, K}(\beta)&=&\frac{1}{N+K}\mathbb{E}_\xi \log Z_{N, K}(\beta),\\
f_{N, K}(\beta)&=&-\frac{1}{\beta}A_{N, K}(\beta).
\end{eqnarray}
We can define also the Boltzmann state for a generic function of the spin variables $F(\sigma, \tau)$:
$$
\omega_{N, K}(F)=Z^{-1}_{N, K}(\beta)\sum_{\sigma,\tau}\left[F(\sigma, \tau)\exp\left(-\beta H_{N,K}(\sigma, \tau)\right)\right].
$$

\ni The main goal of the theory is the control of the free energy in the thermodynamic limit, \textit{i.e.} for $N,K\to\infty$,  when the size of the two parties grows to infinity, such that $N/(N+K)\to \alpha\in(0, 1)$ and $K/(N+K)\to(1-\alpha)\in(0, 1)$. 

\ni We adopt this latter definition of
thermodynamic limit, and thus the thermodynamic functions depend
by the additional parameter $\alpha$, ruling the relative ratio
between the parties:
$$
\lim_{N,K} A_{N, K}(\beta)=A(\alpha, \beta).
$$

\ni At the moment no rigorous proof of the existence of such a limit is known.

\section{Results}

We claim that the free energy of bipartite spin glass models is universal with respect to the choice of the statistical distribution of the quenched noise.

\ni In order to give a rigorous proof of this assertion, we need to recall the following result \cite{CH}:

\begin{lemma}
Let $\xi$ a real r.v. such that, if $g$ is a unit centered gaussian r.v., it is $\mathbb{E}_g g^k=\mathbb{E}_\xi \xi^k$, $\forall k=1...m$, and $\mathbb{E}[|\xi^{m+1}|]$ is finite. Furthermore be $f$ a real function in $C^{m}$, such that $\|f^{(m)}\|_{\infty}=\sup |f^{(m)}|<\infty$. Then
\begin{equation}\label{eq:disf}
|\mathbb{E}[\xi f(\xi)]-\mathbb{E}[\xi^2] \mathbb{E}[f'(\xi)]|\leq \left(\frac{m+1}{m!}\right)\mathbb{E}[|\xi|^{m+1}] \|f^{(m)}\|_{\infty}.
\end{equation}
\end{lemma}

\ni This lemma gives the error we make when we use the gaussian derivative rule for functions of random variables close to be gaussian up to order $m$. Since we want to use Guerra's interpolation technique, it is crucial to have such a lemma that compares derivatives. For example, in the Lindeberg approach, this is replaced by an integral analogue (see Theorem 1.1 in \cite{ch}). The two methods, even though the Lindeberg approach gives a slightly weaker condition on the random interaction, seem to be equivalent: they give similar estimates, and they both depend in a crucial way on the behaviour of the derivatives of the Boltzmann mean value of the spin part of the internal energy (that is, $\omega(\s_i\t_\mu)$, or, in the notation of \cite{ch}, $\omega(\lambda_2)$).

\ni If we take a random matrix $\xi_{ij}$, each entry with the hypothesys of the previous lemma, we define the pressure of the $\xi$-noise model as
\be\label{eq:AxiN}
A_{N+K}^{\xi}=\frac{1}{N+K}\log\sum_{\sigma,\tau} \exp\left( \beta\sqrt{\frac{2}{N+K}}\sum_{i,j=1}^{N,K} \xi_{ij}\sigma_i\tau_j \right),
\ee
while we set
\be\label{eq:AgN}
A_{N+K}^{g}=\frac{1}{N+K}\log\sum_{\sigma,\tau} \exp\left( \beta\sqrt{\frac{2}{N+K}}\sum_{ij=1}^{N,K} g_{ij}\sigma_i\tau_j \right),
\ee
with $g_{ij}$ normal distribuited. So we can prove our first result:
\begin{teorema}
Let $\xi$ be a real r.v., as in the hypothesis of the previous Lemma, with $m\geq2$, and $g$ be a unit centered gaussian. Then we have
\begin{equation}
\left |\E_{g,\xi}\left[A^\xi_{N+K}(\beta)-A^g_{N+K}(\beta)\right]\right|\leq\frac{(\sqrt{2}\beta)^{m+1}\alpha(1-\alpha)}{(N+K)^{(m-1)/2}}\left (m+1\right)\mathbb{E}[|\xi|^{m+1}].
\end{equation}
\end{teorema}
\noindent That is completely analogous with the achievement obtained in \cite{CH}\cite{ch} for the Sherrington-Kirkpatrick model.

\ni Now we have to evaluate fluctuations. We will use an argument based on the cavity technique in order to state the following
\begin{teorema}
Be $\xi_{ij}$ the entries of a random matrix such that, for a fixed $m>0$, $\E[|\xi|^p]$ is bounded for every positive real number $p\leq m$. Then
\begin{equation}
\E\left[|A^{\xi}_{N,K}-A_{N,K}|^p\right]\leq C_p\left(\frac{2\beta^2\alpha(1-\alpha)}{(N+K)}\right)^{p/2}\E[|\xi|^p],
\end{equation}
with $C_p$ an universal constant, depending only by $p$.
\end{teorema}
\ni The combination of the two theorems gives immediately the following
\begin{corollario}
In the hypothesis of the Lemma, for a fixed $m\geq 2$ and a positive $p\leq m+1$ we have
\bea
\E_\xi\left[|A^{\xi}_{N,K}-\E_g A^g_{N,K}|^p\right]&\simeq& \OO{N^{p/2}},\\
\E_g\left[|A^{g}_{N,K}-\E_\xi A^\xi_{N,K}|^p\right]&\simeq& \OO{N^{p/2}}.
\eea
\end{corollario}

\section{Proofs}
\begin{proof}[Proof of the Lemma]
If we expand in Taylor series both the function $f$ and its derivative $f^{(1)}$, we get
\begin{eqnarray}
\left| f(\xi)-f(0)-\xi f^{(1)}(0) - ... - \frac{\xi^{m-1}}{(m-1)!}f^{(m-1)}(0) \right|&\leq&\frac{\xi^{m}}{m!} \|f^{(m)}\|_{\infty}\nn\\
\left| f^{(1)}(\xi)-f^{(1)}(0)-\xi f^{(2)}(0) - ... - \frac{\xi^{m-2}}{(m-2)!}f^{(m-1)}(0) \right|&\leq&\frac{\xi^{m-1}}{(m-1)!} \|f^{(m)}\|_{\infty}\nn.
\end{eqnarray}
Now, since the mean value of $\xi$ is zero, we have
\be\label{eq:lemma-passo1}
\left|\mathbb{E}[\xi f(\xi)]-\mathbb{E}[\xi^2] \mathbb{E}[f'(\xi)]\right|=\left| \mathbb{E}[\xi (f(\xi)-f(0))]-\mathbb{E}[\xi^2] \mathbb{E}[f'(\xi)] \right|.
\ee
Due to the given relation between the $m$ gaussian moments and the ones of $\xi$, we can add and subtract analogous terms in (\ref{eq:lemma-passo1}), according to the Taylor expansions of $f$ and $f^{(1)}$:
\bea
\left|\mathbb{E}[\xi f(\xi)]-\mathbb{E}[\xi^2] \mathbb{E}[f'(\xi)]\right| &=& \left|\mathbb{E}[\xi (f(\xi)-\sum_{i=0}^{m-1}\frac{\xi^i}{i!}f^{(i)}(0)]\right.\nn\\
&-&\left.\mathbb{E}[\xi^2] \mathbb{E}[f'(\xi)-\sum_{i=1}^{m-1}\frac{\xi^{i-1}}{(i-1)!}f^{(i)}(0)]\right| \nn\\
&\leq&\mathbb{E}\left[|\xi| \left|f(\xi)-\sum_{i=0}^{m-1}\frac{\xi^i}{i!}f^{(i)}(0)\right|\right]\nn\\
&+&\mathbb{E}[\xi^2] \mathbb{E}\left[\left|f'(\xi)-\sum_{i=1}^{m-1}\frac{\xi^{i-1}}{(i-1)!}f^{(i)}(0)\right|\right]\nn\\
&\leq&\frac{\mathbb{E}[|\xi|^{m+1}]}{m!}\|f^{(m)}\|_{\infty}+\mathbb{E}[\xi^2]\frac{\mathbb{E}[|\xi|^{m-1}]}{(m-1)!}\|f^{(m)}\|_{\infty}\nn\\
&\leq&\frac{\|f^{(m)}\|_{\infty}}{m!}\left( \mathbb{E}[|\xi|^{m+1}]+m\mathbb{E}[\xi^2] \mathbb{E}[|\xi|^{m-1}] \right)\nn\\
&\leq& \left(\frac{m+1}{m!}\right) \mathbb{E}[|\xi|^{m+1}] \|f^{(m)}\|_{\infty}.\nn
\eea
\end{proof}

\begin{proof}[Proof of Theorem 1]
It is useful to introduce the interpolating partition function:
$$
Z_{N+K}(t)=\E_{\sigma, \tau}\exp\left( \beta\sqrt{\frac{2t}{N+K}}\sum_{ij}\xi_{ij}\sigma_i\tau_j+\beta\sqrt{\frac{2(1-t)}{N+K}}\sum_{ij}g_{ij}\sigma_i\tau_j \right),
$$
and the interpolating pressure
$$
A_{N+K}(t)=\E_{\xi,g}\log Z_{N+K}(t).
$$
It is easily seen that $A_{N+K}(0)=\E_g A_{N+K}^g$ and $A_{N+K}(1)=\E_\xi A_{N+K}^{\xi}$. Furthermore, in virtue of the previous lemma, we have
\be
\frac{d}{dt} A_{N+K}(t)=\frac{\beta}{\sqrt{2(N+K)^{3}}}\left(\sum_{ij}\frac{1}{\sqrt{t}}\E[\xi_{ij}\omega_t(\sigma_i\tau_j)]-\sum_{ij}\frac{1}{\sqrt{1-t}}\E[g_{ij}\omega_t(\sigma_i\tau_j)_t]\right),
\ee
and since
$$
\frac{1}{\sqrt{1-t}}\sum_{ij}\E[g_{ij}\omega_{t}(\sigma_i\tau_j)]=\frac{1}{\sqrt{1-t}}\sum_{ij}\E[\partial_{g}\omega_{t}(\sigma_i\tau_j)]=\frac{1}{\sqrt{t}}\sum_{ij}\E[\partial_{\xi}\omega_{t}(\sigma_i\tau_j)],
$$
bearing in mind the Lemma, we get
$$
\left|\frac{1}{\sqrt{t}}\sum_{ij}\E[\xi_{ij}\omega_{t}(\sigma_i\tau_j)]-\frac{1}{\sqrt{t}}\sum_{ij}\E[\partial_{\xi}\omega_{t}(\sigma_i\tau_j)]\right|\leq\frac{NK}{\sqrt{t}} \left(\frac{m+1}{m!}\right)\mathbb{E}[|\xi|^{m+1}]\E[\|\omega^{m}\|].
$$
Then we easily have
$$
\left|\frac{d}{dt}A_{N+K}(t)\right|\leq\sqrt{N+K}\frac{\beta\alpha(1-\alpha)}{\sqrt{t}}\left(\frac{m+1}{m!}\right)\mathbb{E}[|\xi|^{m+1}]\E[\|\omega^{(m)}\|].
$$
Now we have to estimate the $m$-th derivative of the state $\omega_{t}(\sigma_i\tau_j)$ with respect to $\xi$. If we name $P^{m}$ the cumulant polynomial of degree $m$ in $(\sigma_i\tau_j)$ with respect to the Gibbs measure (for example $P^{3}=\omega((\sigma_i\tau_j)^3)-3\omega((\sigma_i\tau_j)^2)\omega(\sigma_i\tau_j)+2\omega^3(\sigma_i\tau_j)$, see for instance \cite{AC}: the coeffincients are infact the same of Aizenmann Contucci relations in the SK model), it is easy to check that 
\be\label{eq:stima-derivata1}
\omega^{(m)}=\left(\beta\sqrt{2t}\right)^m(N+K)^{-m/2}P^{m+1}
\ee
hence (for dichotomic spin variables)
\be\label{eq:stima-derivata2}
\|\omega^{(m)}(\sigma_i\tau_j)\|\leq\left(\beta\sqrt{2t}\right)^m(N+K)^{-m/2}(m+1)!.
\ee
Therefore, it follows that
\bea
\left|\E\left[A^{\xi}_{N+K}-A^{g}_{N+K}\right]\right|&\leq&\int_0^1dt \left|\frac{d}{dt} A_{N+K}(t)\right|\nn\\
&\leq& \frac{(\sqrt{2}\beta)^{m+1}\alpha(1-\alpha)}{(N+K)^{(m-1)/2}}\int_0^1dt \frac{t^{(m-1)/2}}{2} \mathbb{E}[|\xi|^{m+1}] \nn\\
&=&\frac{(\sqrt{2}\beta)^{m+1}\alpha(1-\alpha)}{(N+K)^{(m-1)/2}}(m+1)\mathbb{E}[|\xi|^{m+1}]\nn,
\eea
and the theorem is proven.
\end{proof}

\begin{proof}[Proof of Theorem 2]
In primis we construct an increasing sequence of subsets in the set of spin variables. Given a $NK$-size system, we cover it with sets labeled by $h$, such that for each $h$ there is a subset of size $N_hK_h<NK$,  and when $h$ increases by one, we have alternatively either $N_{h+1}=N_{h}+1$, $K_{h+1}=K_{h}$, or $N_{h+1}=N_h$, $K_{h+1}=K_h+1$, $h=0, ..., NK$. Furthermore we have $N_0=K_0=0$ and $N_{NK}=N$, $K_{NK}=K$. This sequence induces a natural filtration on the $\xi$ variables, $F_h=\{\xi_1,...,\xi_h\}$, \textit{i.e.} the noise occurring to define a $N_hK_h$-size bipartite system. 

\ni Then we proceed with a cavity argument. We can write the hamiltonian of the system as the hamiltonian of a smaller one, with suitable additional terms, namely the cavity field. This is a standard approach in spin glasses \cite{MPV}\cite{talabook}. Thus, bearing in mind our sequence, we write
\bea
-\beta H_{N,K}&=&-\beta_h H_{N_h,K_h}+\beta_h\sqrt{2(1-\alpha_h)}\sum_{i=N_h+1}^N\bar\sigma_i h^i(\tau)\nn\\
&+&\beta_h\sqrt{2\alpha_h}\sum_{j=K_h+1}^K\bar\tau_j h^j(\sigma)+\beta\sqrt{\frac{2}{N+K}}\sum_{i=N_h+1}^N\sum_{j=K_h+1}^K\bar\xi_{ij}\bar\sigma_i\bar\tau_j,\nn
\eea
where $\beta_h=\beta\sqrt{(N_h+K_h)/(N+K)}$, $\alpha_h=N_h/(N_h+K_h)$, as usual $h^i(\tau)=K_h^{-1/2}\sum_{j=1}^{K_h}\xi_{ij}\tau_j$ and $h^j(\sigma)=N_h^{-1/2}\sum_{i=1}^{N_h}\xi_{ij}\sigma_i$ are the cavity fields, and we indicate with a bar the spin and the noise variable outside of the $h$-th sequence and filtration.
Hence we have
\bea
\log Z_{N,K}(\beta)&=&\log Z_{N_h,K_h}(\beta_h)\nn\\
&+&\log \omega_{N_h, K_h}\left(\E_{\bar\s_{N_h+1}...\bar\s_N}\E_{\bar\t_{K_h+1}...\bar\t_K}\exp\left[ \beta_h\sqrt{2(1-\alpha_h)}\sum_{i=N_h+1}^N\bar\sigma_i h^i(\tau)\right.\right.\nn\\
&+&\left.\left.\beta_h\sqrt{2\alpha_h}\sum_{j=K_h+1}^K\bar\tau_j h^j(\sigma)+\beta\sqrt{\frac{2}{N+K}}\sum_{i=N_h+1}^N\sum_{j=K_h+1}^K\bar\xi_{ij}\bar\sigma_i\bar\tau_j \right]\right).\nn
\eea

\ni Now we can define the function
\bea
\psi_h(\bar\xi)&=&\E_{\xi_1...\xi_h} \log \omega_{N_h, K_h}\left(\E_{\bar\s_{N_h+1}...\bar\s_N}\E_{\bar\t_{K_h+1}...\bar\t_K}\exp\left[ \beta_h\sqrt{2(1-\alpha_h)}\sum_{i=N_h+1}^N\bar\sigma_i h^i(\tau)\right.\right.\nn\\
&+&\left.\left.\beta_h\sqrt{2\alpha_h}\sum_{j=K_h+1}^K\bar\tau_j h^j(\sigma)+\beta\sqrt{\frac{2}{N+K}}\sum_{i=N_h+1}^N\sum_{j=K_h+1}^K\bar\xi_{ij}\bar\sigma_i\bar\tau_j \right]\right),\nn
\eea
that is a suitable generalization of the well known cavity function introduced in spin glass theory many years ago \cite{Glocarno}. It is worthwhile to notice that $\psi_0=\log Z_{N,K}$, and $\psi_{NK}=0$. So we have that
$$
\E [\log Z_{N,K}|F_h]=(N_h+K_h)A_{N_h,K_h}+\psi_h(\bar \xi_{h}).
$$
Now, following \cite{CH}, we can introduce the martingale
\be\label{eq:deltah}
\Delta_h=\E [\log Z_{N,K}|F_h]-\E [\log Z_{N,K}|F_{h-1}],
\ee
and notice that
\be\label{eq:deltah-dec}
\left|\frac{1}{N+K}\sum_{h=1}^{NK}\Delta_h\right|=|A^{\xi}_{N,K}-A_{N,K}|.
\ee
At this point we need a bound on $\Delta_h$:
\bea
|\Delta_h|&=&|\psi_h-\psi_{h-1}+\log Z_{N_h, K_h}-\log Z_{N_{h-1},K_{h-1}}|\nn\\
&=&\log \omega^*\left( e^{\beta\sqrt{\frac{2}{N+K}}\bar\xi_h \bar\sigma \bar\tau} \right)
\leq \beta\sqrt{\frac{2}{N+K}}|\xi|\label{eq:stima-delta},
\eea
hence
\be\label{eq:Delta^2}
\Delta_h^2\leq \frac{2\beta^2}{N+K}|\xi|^2.
\ee
Here $\omega^*$ denotes the expectation value with respect to an ausiliar Gibbs measure, with weights 
\bea
& &\exp\left(-\beta_h H_{N_{h-1},K_{h-1}}+\beta_h\sqrt{2(1-\alpha_h)}\sum_{i=N_h+1}^N\bar\sigma_i h^i(\tau)\right.\nn\\
&+&\left.\beta_h\sqrt{2\alpha_h}\sum_{j=K_h+1}^K\bar\tau_j h^j(\sigma)+\beta\sqrt{\frac{2}{N+K}}\sum_{i=N_h+1}^N\sum_{j=K_h+1}^K\bar\xi_{ij}\bar\sigma_i\bar\tau_j \right).\nn
\eea
Now we use a martingale moments inequality (see for instance \cite{IW}), in order to state
\be\label{eq:mart}
\E\left|\frac{1}{N+K}\sum_{h=1}^{NK}\Delta_h\right|^p\leq C_p\E\left|\frac{1}{(N+K)^2}\sum_{h=1}^{NK}\Delta^2_h\right|^{p/2},
\ee
with $C_p$ an universal constant independent on $\xi$. Therefore, putting together (\ref{eq:deltah-dec}), (\ref{eq:Delta^2}) and (\ref{eq:mart}):
$$
\E\left[|A^{\xi}_{N,K}-A_{N,K}|^p\right]\leq C_p\left(\frac{2\beta^2\alpha(1-\alpha)}{(N+K)}\right)^{p/2}\E[|\xi|^p].
$$
\end{proof}

\section{Further Remarks and Open Problems}

In this work we have established the invariance of the free energy of bipartite spin glass, with respect to the choice of the statistical distribution of the random interaction in thermodynamic limit. Our method is on the same line of \cite{CH}, and anyway both are based on Guerra's interpolation and cavity field approach. Our results can be interpreted as follows: at first we have shown that, under certain hypothesis on the random interaction, the free energy converges in distribution to the one with gaussian interaction (Theorem 1); then we have shown that it converges also in $L_p$, for suitable values of $p$ (Theorem 2). Of course all that holds provided the proof of the existence of the thermodynamic limit of the free energy for bipartite models of spin glasses, that remains an open problem.

\ni We can make some other consideration. At first we notice that in bipartite models with $\pm 1$ spin, we can improve estimate (\ref{eq:stima-delta}), in order to get faster rate of convergence, but for less values of $p$:
\be\label{eq:stima-migl}
\log \omega^*\left( e^{\beta\sqrt{\frac{2}{N+K}}\bar\xi_h \bar\sigma \bar\tau} \right)
\leq \beta^2\frac{2}{N+K}\xi^2
\ee
This is convenient in the models of interest, where the random interaction has all moments bounded (\textit{e.g.} $\pm 1$, gaussian). So we can give an other version of Theorem 2:
\begin{t*}[\textbf{2 v2}]
If we deal with dichotomic spin variables, in the hypothesis of Theorem 2 but for $0<p<(m+1)/2$, we have
$$
\E\left[|A^{\xi}_{N,K}-A_{N,K}|^p\right]\leq C_p\left(\frac{2\beta^2\sqrt{\alpha(1-\alpha)}}{(N+K)}\right)^{p}\E[|\xi|^{2p}],
$$
with $C_p$ an universal constant, depending only by $p$.
\end{t*}

\ni The last theorem holds also if only one of the party is made by $\pm 1$ spin, and the other is formed by soft bounded spin. Anyway we have the following

\begin{corollario}
Consider a model of bipartite spin glass with spin $\pm 1$. Provided the thermodinamic limit for the free energy exists, we have that $\forall p>0$ $A_{N,K}^g\overset{L_p}{\longrightarrow} A(\alpha, \beta)$ (gaussian random interaction) and $A_{N,K}^{\pm}\overset{L_p}{\longrightarrow} A(\alpha, \beta)$ (dichotomic random interaction). Furthermore, in both cases and $\forall p>0$, the convergence rate is $\max\left(\OO{(N+K)^p},|A_{N,K}(\beta)-A(\alpha,\beta)|^p\right)$.
\end{corollario}

\ni For unbounded spin the model is the most interesting: if we take a bipartite model with $N$ $\pm 1$ and $K$ gaussian spin, interacting via the hamiltonian
$$
H_N=-\frac{1}{\sqrt{N}}\sum_{i\mu} \xi_{i\mu}\s_i\t_\mu,
$$
this turns out to be equivalent to the Hopfield Model for neural network \cite{talabook}\cite{Bovbook}. We notice that the normalization is different to the one used in this paper, hence now $\alpha=K/N$. Originally the model was introduced with $\pm 1$ random interaction. Later also its gaussian version has been studied \cite{hop}\cite{BGG1}.

\ni Thus we would like to extend in part our results to this model; unfurtunately this seems very hard to do. The main difficulty is in Theorem 1: in both the approaches for proving universality in spin glasses, the one based on the Wick rule, or interpolation method, (\cite{talbern}, \cite{CH}, ours), or the one based on Lindeberg argument \cite{ch}, it is crucial to get an estimate on the derivatives of $\omega (\s_i\t_\mu)$; that is trivial in our model, or even in SK model, but it has a deep significance in the Hopfield Model. Actually we know that universality (in a strict sense) in the Hopfield Model does not hold: it suffices to match the results one can get for the model with gaussian or dichotomic interaction \cite{AGS}\cite{BGG1}. Infact in general we have 
$$
A^g_{N,K}(\beta)\geq A_{N,K}\left(\frac{2\beta}{\pi}\right),
$$
and furthermore the models are of course different in the limit $\alpha\to0$ \cite{pr}.

\ni So a sharp estimate on the derivatives of $\omega (\s_i\t_\mu)$ would give a region where universality is expected. This is anyway an open problem.

\ni On the other hand, we can easily extend Theorem 2 (better in its second version, performing an integration over gaussian variables in order to get (\ref{eq:stima-migl})) to the Hopfield Model, with both gaussian or dichotomic interaction:
\begin{t*}[\textbf{2} for the Hopfield Model]
In the hypothesis of Theorem 2 but for $0<p<(m+1)/2$, we have
$$
\E\left[|A^{\xi}_{N,K}-A_{N,K}|^p\right]\leq C_p\left(\frac{\beta\sqrt{\alpha}}{N}\right)^{p}\E[|\xi|^{2p}],
$$
with $C_p$ an universal constant, depending only by $p$.
\end{t*} 
\ni and we get immediately the corollary
\begin{corollario}
Provided the thermodinamic limit for the free energy of the Hopfield Model exists, for both the models with gaussian or $\pm 1$ random interaction, we have that $\forall p>0$ $A_{N,K}^g\overset{L_p}{\longrightarrow} A^{g}(\alpha, \beta)$ (gaussian random interaction) and $A_{N,K}^{\pm}\overset{L_p}{\longrightarrow} A^{\pm}(\alpha, \beta)$ (dichotomic random interaction). Furthermore, in both cases and $\forall p>0$, the convergence rate is $\max\left(\OO{N^p},|A_{N,K}(\beta)-A(\alpha,\beta)|^p\right)$.
\end{corollario}

\ni Of course the proof of the existence of the thermodynamic limit for the free energy of the Hopfield Model is as hard as in the bipartite spin glass. But for the original Hopfield Model we can do even better, proving convergence in the $\sup$ norm. Infact by (\ref{eq:deltah}) and (\ref{eq:stima-migl}) we are sure that the difference between the pressure and its expectation value is bounded uniformly in $\xi$:
$$
\sup_{\xi=\pm1}|A^{\xi}_{N,K}(\xi)-A_{N,K}|\leq\sup_{\xi=\pm1}|\xi|^2\alpha\beta^2=\alpha\beta^2.
$$
Furthermore
$$
\|A^{\xi}_{N,K}-A_{N,K}\|_p\leq \sqrt[p]{C_p}\frac{\alpha\beta^2}{N}
$$
so
\bea
\lim_p\|A^{\xi}_{N,K}-A_{N,K}\|_p&=&\|A^{\xi}_{N,K}-A_{N,K}\|_{\infty}\nn\\
&=&\lim_p\sqrt[p]{C_p}\frac{\alpha\beta^2}{N}\nn\\
&\leq&\frac{\alpha\beta^2}{N}C
\eea
since, for $p>1$, $C_p\leq C p$ for a certain costant $C>0$ \cite{IW}. Thus as a direct conseguence of our work we have the following final
\begin{teorema}
Provided the existence of the thermodynamic limit for the Hopfield Model, we have that the pressure converges in $L_{\infty}$ norm to its mean value, with rate
$$
\max\left(\frac{1}{N}, |A_{N,K}(\beta)-A(\alpha,\beta)|\right).
$$
\end{teorema}

\ni This result is in the same wake of the ones in \cite{Bovbook}\cite{talabook}. 

\bigskip
\noindent
{\bf Acknowledgements}
\newline
I am grateful to Francesco Guerra for his precious scientific guide. Furthermore, I warmly thank Adriano Barra, Renato Luc\`a, Daniele Tantari for many fruitful discussions, and Chiara Saffirio for reference \cite{IW}.

\end{document}